\documentclass{amsart}
\usepackage[latin1]{inputenc}
\usepackage{amssymb,amsmath}

\newenvironment{dem}[1]%
    {\par\noindent{\scshape Proof:\ }#1}%
    {\mbox{}\hfill$\square$\par\bigskip\par}
\newtheorem{teo}{Theorem}[section]
\newtheorem{defi}[teo]{Definition}
\newtheorem{prop}[teo]{Proposition}
\newtheorem{coro}[teo]{Corollary}
\newtheorem{lema}[teo]{Lemma}
\newcommand{\To}{\longrightarrow}

\begin{document}
\title{Boolean Metric Spaces and Boolean Algebraic Varieties}

\author{Antonio Avilés}\thanks{Author supported by FPU grant of SEEU-MECD, Spain.}
\email{avileslo@um.es}
\address{Departamento de Matemáticas. Universidad de
Murcia, 30100 Murcia (Spain)}

\begin{abstract}The concepts of Boolean metric space and convex
combination are used to characterize polynomial maps $A^{n}\To
A^{m}$ in a class of commutative Von Neumann regular rings
including $p$-rings, that we have called CFG-rings. In those
rings, the study of the category of algebraic varieties (i.e. sets
of solutions to a finite number of polynomial equations with
polynomial maps as morphisms) is equivalent to the study of a
class of Boolean metric spaces, that we call here CFG-spaces.\end{abstract}

\maketitle

\subsection*{Notations and conventions}

Throughout this work, $(B,+,\cdot)$ will be a Boolean ring where
the operation $a\vee b = a+b+ab$ is the analogue for set union,
the order $a\leq b \Leftrightarrow ab=a$ is the analogue for set
inclusion and for each $a\in B$, $\bar{a}=a+1$ is the analogue for
the set complement of $a$.\\

All rings will be commutative with
identity. Regular ring will mean here commutative Von Neumann
regular ring, i.e. a (commutative) ring for which any principal
ideal is generated by an idempotent, also known as absolutely flat rings, see~\cite{Goodearl}, \cite{Vra1}.
Unless otherwise stated, $A$ will
be a regular ring, $B(A)$
will denote the set of the idempotent elements of $A$ and
$e:A\longrightarrow B(A)$ will be the map that sends each $a\in A$
to the only idempotent $e(a)\in B(A)$ such that $aA=e(a)A$. The
set $B(A)$ has a structure of Boolean ring with product inherited
from $A$ and with the sum $a\tilde{+}b=(a-b)^{2}$. For
$a_{1},\ldots,a_{n}\in B(A)$ with $a_{i}a_{j}=0$ for $i\neq j$, it
holds
$a=a_{1}+\cdots+a_{n}=a_{1}\tilde{+}\cdots\tilde{+}a_{n}=a_{1}\vee\cdots\vee
a_{n}$. In this case we will denote $a$ by
$a_{1}\oplus\cdots\oplus a_{n}$.\\

Given a prime $p\in\mathbf{Z}$, a $p$-ring is a ring $A$ for which
$px=0$ and $x^{p}=x$ for all $x\in A$. In particular, a Boolean
ring is a 2-ring. Any $p$-ring is a regular ring with
$e(x)=x^{p-1}$.\\

An algebraic variety over a ring $A$ is a set $U\subset A^n$ which is 
the set of solutions to a finite number of polynomial equations. If 
$U\subset A^n$ and $V\subset A^m$ are algebraic varieties, a map $f:U\To 
V$ is called a polynomial map if there are polynomials $f_1,\ldots,f_m\in A[X_1,\ldots 
X_n]$ such that $f(x) = (f_1(x),\ldots,f_m(x))$. When $A=B$ is a Boolean 
ring the usual terms are Boolean domain and Boolean transformation, see~\cite{Rud} and~\cite{Ru2}.

\subsection*{Introduction}

Boolean metric spaces (Definition~\ref{espacio metrico}) appeared
in several works in the 1950's and 1960's \cite{Blu1}, \cite{Blu2},
\cite{Ell1}, \cite{Ell2}, \cite{Har} and  \cite{Mel1}, where some authors investigated the
analogue for some topics in Geometry such as \emph{betweeness},
motions or topology in some of those spaces, as
Boolean algebras and some rings where a suitable Boolean metric
could be defined. In some papers, \cite{Bat}, \cite{Mel2}, \cite{Mr2} and
\cite{Zem}, a special attention was paid to $p$-rings, that admit
a metric space structure over its ring of idempotents. In fact, if $A$ is 
a regular ring, then $A^n$ is a Boolean metric space over $B(A)$ with the 
distance $$d((x_{1},\ldots,x_{n}),(y_{1},\ldots,y_{n})) =
e(x_{1}-y_{1})\vee\cdots\vee e(x_{n}-y_{n}).$$

We
will show the close relation that exists between the theory of
Boolean metric spaces and the Algebraic Geometry over CFG-rings. We define a regular ring $A$ to 
be an CFG-ring if there are $x_1,\ldots,x_n$ in $A$ such that any element in $A$ 
is of the form $\sum a_{i}x_{i}$ where $a_1,\ldots a_n\in B(A)$ and 
$a_1\oplus\cdots\oplus a_n=1$.\\

In sections 1 and 2 we develop some tools concerning the 
structure of Boolean metric spaces, while in sections 3 and 4 the 
main results are exposed. Namely, in section 3 we prove that if $A$ is a CFG-ring and 
$U$ is a subset of $A^n$, then $U$ is an algebraic variety if and only if 
there are $x_1,\ldots,x_n$ in $U$ such that any element of $x\in U$ is of 
the form $x=\sum_1^n a_i x_i$ where $a_1,\ldots,a_n\in B(A)$ and 
$a_1\oplus\cdots\oplus a_n=1$, if and only if there is distance-preserving 
bijection from $U$ onto an algebraic variety $V\subset A^m$. Also, if 
$U\subset A^n$ and $V\subset A^m$ are algebraic varieties over $A$ and 
$f:U\To V$ is a map, the following are equivalent:\begin{enumerate}
\item $f$ is a polynomial map.
\item $d(f(x),f(y))\leq d(x,y)$ for all $x,y$ in $U$.
\item $f(\sum_1^n a_i x_i) = \sum_1^n a_i f(x_i)$ for all $x_1,\ldots,x_n$ in 
$U$ and for all $a_1,\ldots,a_n$ in $B(A)$ with $a_1\oplus\cdots\oplus a_n=1$.
\end{enumerate}

Thus, the category of algebraic varieties over an CFG-ring is equivalent 
to the category of those Boolean metric spaces over $B(A)$ that are 
isometric to some algebraic variety, that we have called CFG-spaces. Some special
cases of these implications were known for $A=B$ a Boolean ring: that 1 is 
equivalent to 3 when $U=B^n$, $V=B$ is in Theorem 4.2 in~\cite{Rud}, and that 1 is equivalent 
to 2 when $U=V=B$ was observed in~\cite{Mr1}.\\

 In section 4, we present 
a classification of the Boolean metric spaces over a Boolean ring $B$, 
which is a classification of the algebraic varieties over a CFG-ring. We 
associate to each of those spaces a finite decreasing sequence of nonzero 
elements of $B$ such that two spaces are isometric if and only if 
they have the same associated sequence.\\

The author wishes to thank professors Juan Martínez and Manuel
Saorín, from University of Murcia, and Sergiu Rudeanu, from University of Bucharest, for their support and stimulus, and for their help in the redaction of this article.

\section{Boolean metric spaces}

\subsection{Basic definitions and examples}

\begin{defi}\label{espacio metrico}
Let $X$ be a set. A map $d:X\times X\longrightarrow B$ is said to
be a Boolean metric if the following axioms hold, for all
$x,y,z\in X$:
\begin{enumerate}
\item $d(x,y)=0$ if and only if $x=y$.
\item $d(x,y) = d(y,x)$.
\item $d(x,z) \leq d(x,y)\vee d(y,z)$.
\end{enumerate}
In that case, we will say that $(X,d)$ is a metric space over $B$.
\end{defi}

In the above definition, axiom 3 can be substituted by any of the
following:
\begin{itemize}
\item[3'.] $d(x,z)\overline{d(y,z)}\leq d(x,y)$
\item[3''.] $d(x,z) + d(z,y)\leq d(x,y)$
\end{itemize}

Some suitable subsets of modules possess structure of Boolean metric 
space. We have called these subsets metrizable. Recall that the 
annihilator of an element $x$ of a module over the ring $A$ is the ideal 
$Ann(x) = \{a\in A : ax=0\}$.

\begin{defi}
Let $A$ be a regular ring, $M$ a module over $A$ and $X$ a subset of $A$. 
The set $X$ will be said to be metrizable if for each $x,y\in X$ the ideal 
$Ann(x-y)$ is a principal ideal of $A$.
\end{defi}

If $X$ is a metrizable subset of $M$, for each $x,y\in M$ the ideal 
$Ann(x-y)$ has a unique idempotent generator, say $a_{xy}\in B(A)$. Then, the map $d(x,y) = \overline{a_{xy}}$ is a 
Boolean metric on $X$, called the modular metric on $X$. 
Triangular inequality follows from 
$$Ann(x-y)\cap Ann(y-z)\subseteq Ann(x-z).$$

For every $a\in A$ we have $Ann(a)=\overline{e(a)}A$, so $A$ is a
metrizable subset of itself and its modular metric is given by
$d(x,y)=e(x-y)$. This is the same metric on $A$ as defined in~\cite{Mel1}.\\

Furthermore, for every $n\in\mathbf{N}$, $A^{n}$
is also a metrizable subset of itself and its modular metric is
given by $$d((x_{1},\ldots,x_{n}),(y_{1},\ldots,y_{n})) =
e(x_{1}-y_{1})\vee\cdots\vee e(x_{n}-y_{n}).$$ This is a particular case 
of the following general construction:

\begin{defi}\label{ejemplo X=producto}
Let $(X_{1},d_{1}),\ldots,(X_{n},d_{n})$ be metric spaces over
$B$. Then $(X_{1}\times\cdots\times X_{n},d)$ is also a metric
space over $B$ with $$d((x_{1},\ldots,x_{n}),(y_{1},\ldots,y_{n}))
:= d_1(x_{1},y_{1})\vee\cdots\vee d_n(x_{n},y_{n})$$ This space will
be called the product space of the spaces $(X_{i},d_{i})$ and $d$
will be called the product metric of the metrics $d_{i}$.
\end{defi}

The formation of products is compatible with modular metrics:

\begin{prop}\label{producto de medibles es medible}
Let $S_{i}$ be a metrizable subset of the $A$-module $M_{i}$, for
$i=,1\ldots,n$. Then $S = S_{1}\times\cdots\times S_{n}$ is a
metrizable subset of $M_{1}\times\cdots\times M_{n}$ and the
modular metric in $S$ equals the product metric of the modular
metrics in the $S_{i}$'s.
\end{prop}

\begin{dem}
Call $d_{i}$ the modular metric in $S_{i}$. For each
$x=(x_{1},\ldots,x_{n})$ and $y= (y_{1},\ldots,y_{n})$ in $S$
\begin{eqnarray*} Ann(x-y) &=& \bigcap_{i=1}^{n}Ann(x_{i}-y_{i}) =
\bigcap_{i=1}^{n}\overline{d_{i}(x,y)}A\\ &=&
\left(\prod_{i=1}^{n}\overline{d_{i}(x,y)}\right)A =
\left(\overline{\bigvee_{i=1}^{n}d_{i}(x,y)}\right)A.\end{eqnarray*}

\end{dem}

\begin{defi}
Let $X$ and $Y$ be Boolean metric spaces over $B$. A map
$f:X\longrightarrow Y$ is said to be
\begin{enumerate}
\item contractive if $d(f(x),f(y))\leq d(x,y)$ for all $x,y\in
X$.
\item an immersion if $d(f(x),f(y)) = d(x,y)$ for all $x,y\in
X$.
\item an isometry if it is a bijective immersion.
\end{enumerate}
\end{defi}

Contractive maps play the rôle of morphisms in the category of 
Boolean metric spaces over $B$, while isometries are the isomorphisms.
Observe that immersions are always into.\\

\begin{teo}\label{inclusion en un modulo}
Every metric space $X$ over $B$ is isometric to a metrizable
subset of a $B$-module. Furthermore, if we fix $x_{0}\in X$ there
is a metrizable subset $S$ of a $B$-module $M$ such that $0\in S$
and an isometry $g:X\longrightarrow S$ such that $g(x_{0})=0$.
\end{teo}

\begin{dem}
We define $f:X\To B^{X}$ by $f(x) = (d(x,z))_{z\in X}$. To prove
that $f(X)$ is metrizable and that $f:X\To f(X)$ is an isometry,
it is enough to see that $Ann(f(x)-f(y))=\overline{d(x,y)}B$ for
all $x,y\in X$. If $a\in
Ann(f(x)+f(y))$ then, $$a(d(x,z)+d(y,z))_{z\in X}=0,$$ so for $z=x$,
we have $ad(y,x)=0$ and therefore $a\leq \overline{d(x,y)}$.
Conversely, suppose $a\in \overline{d(x,y)}B$, then 
$$a(d(x,z)+d(z,y))\leq ad(x,y) = 0,$$ for
all $z\in X$, so $a\in Ann(f(x)+f(y))$. For the last assertion,
take $h:f(X)\To f(X)+f(x_{0})$ given by $h(x)=x+f(x_{0})$. Then,
$h$ is an isometry between $f(X)$ and the metrizable set
$S=f(X)+f(x_{0})$ because $Ann(h(x)-h(y))=Ann(x-y)$ for all $x,y$.
Hence, the map $g=h\circ f$ is an isometry between $X$ and $S$ that
verifies $g(x_{0})=0$.
\end{dem}

\subsection{Convex combinations and convex closures}

Unless otherwise stated, $X$ will be a metric space over $B$.

\begin{defi}
Let $x_{1},\dots,x_{n}\in X$ and let $a_{1},\dots,a_{n}\in B$ such
that
    $a_{1}\oplus\dots\oplus a_{n}=1$. We will say that $x\in X$ is a convex combination of
    $x_{1},\dots,x_{n}$ with coefficients $a_{1},\dots,a_{n}$ if
    $a_{i}d(x,x_{i})=0$ for $i=1,\ldots,n$.
\end{defi}

\begin{prop}\label{convexidad de la distancia}
    If $x\in X$ is a convex combination of $x_{1},\dots,x_{n}$ with
    coefficients $a_{1},\dots,a_{n}$, then for all $y\in X$
    $$ d(x,y) = \bigoplus_{i=1}^{n}a_{i}d(x_{i},y) $$
\end{prop}

\begin{dem}
For all $i=1,\ldots,n$, since $a_{i}d(x,x_{i})=0$, we have
$a_{i}d(x_{i},y) = a_{i}(d(x,x_{i})+d(x_{i},y)) \leq a_{i}d(x,y)
\leq a_{i}(d(x,x_{i})\vee d(x_{i},y)) = a_{i}d(x_{i},y)$, so
$a_{i}d(x,y)=a_{i}d(x_{i},y)$ and hence, we have $d(x,y)
=(\sum_{i}a_{i})d(x,y) = \sum_{i}a_{i}d(x_{i},y)$.
\end{dem}

\begin{prop}\label{unicidad de las combinaciones convexas}
    If $x$ and $y$ are convex combinations of $x_{1},\dots,x_{n}$
    with coefficients
    $a_{1},\dots,a_{n}$, then $x=y$.
\end{prop}

\begin{dem}
By Proposition \ref{convexidad de la distancia} $$d(x,y)
=
\sum_{i=1}^{n}a_{i}d(x,x_{i})=\sum_{i=1}^{n}a_{i}\sum_{j=1}^{n}a_{j}d(x_{j}x_{i})
= \sum_{i=1}^{n}\sum_{j=1}^{n}a_{i}a_{j}d(x_{j},x_{i})$$ Note that
if $i\neq j$ then $a_{i}a_{j}=0$ and if $i=j$ then
$d(x_{j},x_{i})=0$, so all the terms in the above sum are zero,
and therefore $d(x,y)=0$.
\end{dem}

\begin{lema}\label{envoltura de medible}
Let $S$ be a metrizable subset of an $A$-module $M$. Then,
$$conv(S) = \{a_{1}x_{1}+\cdots+a_{n}x_{n}\in M : x_{i}\in S \
a_{i}\in B(A) \ \bigoplus_{i}a_{i}=1\}$$ is also a metrizable
subset of $M$.
\end{lema}

\begin{dem}
Take $x,y\in conv(S)$, $x=\sum_{1}^{n}a_{i}x_{i}$ and
$y=\sum_{1}^{m}b_{j}y_{j}$. Call $c_{ij}=a_{i}b_{j}$. Then, $\bigoplus_{i,j}c_{ij}=1$ and
$x=\sum_{i,j}c_{ij}x_{i}$ and $y=\sum_{i,j}c_{ij}y_{j}$. Hence,
$Ann(x-y)=Ann(\sum_{i,j}c_{ij}(x_{i}-y_{j}))=
\sum_{i,j}c_{ij}Ann(x_{i}-y_{j})$ which is a principal ideal because every
$Ann(x_{i}-y_{j})$ is principal (recall that, for regular rings,
any finitely generated ideal is principal).
\end{dem}

The following proposition will show that, when $X$ is a metrizable
subset of a module, convex combinations in $(X,d)$ are exactly the
corresponding linear combinations in the module.

\begin{prop}\label{combinacion convexa=combinacion lineal}
Let $S$ be a metrizable subset of an $A$-module,
$x,x_{1},\ldots,x_{n} \in S$ and $a_{1},\ldots,a_{n} \in B$ such
that $\bigoplus_{i=1}^{n}a_{i}=1$. Then, $x$ is a convex
combination of $x_{1},\ldots,x_{n}$ with coefficients
$a_{1},\ldots,a_{n}$ if and only if
$x=a_{1}x_{1}+\cdots+a_{n}x_{n}$.
\end{prop}

\begin{dem}
Suppose $x=a_{1}x_{1}+\cdots+a_{n}x_{n}$. We must check that, for
each $i\in\{1,\ldots,n\}$, $a_{i}d(x,x_{i})=0$. It is clear that
$a_{i}\in Ann(x-x_{i})=\overline{d(x,x_{i})}A$, so $a_{i}\leq
\overline{d(x,x_{i})}$ and hence, $a_{i}d(x,x_{i})=0$.

Conversely, suppose $x\in S$ is a convex combination of
$x_{1},\ldots,x_{n}$ with coefficients $a_{1},\ldots,a_{n}$ . Let
$y=\sum_{1}^{n}a_{i}x_{i}\in conv(S)$, which is metrizable, by
Lemma~\ref{envoltura de medible}. The implication that we have
already proved, tells us that $y$ is a convex combination of
$x_{1},\ldots,x_{n}$ with coefficients $a_{1},\ldots,a_{n}$ in
$conv(S)$. The same holds for $x$, so by Proposition~\ref{unicidad
de las combinaciones convexas}, $x=y$.
\end{dem}

In general, in any metric space $X$, we will denote by
$\sum_{i=1}^{n}a_{i}x_{i}$ or by $a_{1}x_{1}+\cdots+a_{n}x_{n}$
the convex combination of $x_{1},\ldots,x_{n}$ with coefficients
$a_{1},\ldots,a_{n}$, if it exists.\\

Recall that Theorem~\ref{inclusion en un modulo} allows us to identify any
metric space $X$ over $B$ with a metrizable subset of a
$B$-module, and then, by Theorem~\ref{combinacion
convexa=combinacion lineal}, convex combinations are just the
corresponding linear combinations in the module and the metric is
the modular metric.\\

Contractive maps can be
characterized as those that preserve convex combinations.

\begin{teo}\label{convexa sii contractiva}
For a map $f:X\longrightarrow Y$ between two metric spaces $(X,d)$
and $(Y,d')$ the following are equivalent:\begin{enumerate}\item
$f$ is contractive. \item For all $x,x_{1},\ldots,x_{n}\in X$ and
$a_{1},\ldots,a_{n}\in B$ with $\bigoplus_{1}^{n}a_{i}=1$, if
$x=\sum a_{i}x_{i}$, then $f(x)=\sum
a_{i}f(x_{i})$.\end{enumerate}
\end{teo}

\begin{dem}
($1\Rightarrow 2$) Let $x = \sum_{i}a_{i}x_{i}$. Then, for every
$i$, we have $0=a_{i}d(x,x_{i}) \geq a_{i}d(f(x),f(x_{i}))$ , so
$f(x)=\sum a_{i}f(x_{i})$.

  $(2 \Rightarrow 1)$ Given $x,y \in X$, $x
= d(x,y)x + \overline{d(x,y)}y$. Hence, by our assumption
    $f(x) = d(x,y)f(x)+\overline{d(x,y)}f(y)$ and making use of
    Proposition~\ref{convexidad de la distancia} we have finally
    \begin{eqnarray*} d(f(x),f(y)) &=& d(x,y)d(f(x),f(y)) +
    \overline{d(x,y)}d(f(y),f(y))\\
    &=& d(x,y)d(f(x),f(y))\end{eqnarray*} and therefore $d(f(x),f(y))\leq d(x,y)$.
\end{dem}

Given $x_{1},\ldots,x_{n}\in X$ and $a_{1},\ldots,a_{n}\in B$ with
$\bigoplus a_{i}=1$, there may exist no convex combination of the
$x_{i}$'s with coefficients $a_{i}$'s. So we have the next
definition:

\begin{defi}
A metric space $(X,d)$ over $B$ is said to be convex if given any
$x_{1},\ldots,x_{n}\in X$ and any $a_{1},\ldots,a_{n}\in B$ with
$\bigoplus a_{i}=1$, there exists in $X$ the convex combination of
the $x_{i}$'s with coefficients the $a_{i}$'s.
\end{defi}

This notion of convexity is different from the defined in~\cite{Blu2}.

\begin{defi} A convex closure of a metric space $X$ is a convex metric
space $Y\supseteq X$ such that any element in $Y$ is a convex
combination of elements of $X$.\end{defi}

Observe that every metric space $X$ over $B$ has a convex closure, because 
$X$ is isometric to a metrizable subset $S$ of $B$-module and in this 
case, the set $conv(S)$ of
Lemma~\ref{envoltura de medible} is a convex closure of $S$.

\begin{teo}\label{extension a la envoltura convexa}
Let $X\subseteq\bar{X}$ and $Y\subseteq\bar{Y}$ be convex
closures. Each contractive map $f:X\longrightarrow Y$ extends to a
unique contractive map $\bar{f}:\bar{X}\longrightarrow \bar{Y}$.
Furthermore,
\begin{enumerate}\item $\bar{f}$ is immersion if and only if $f$
is, and if $f$ is isometry, so is $\bar{f}$.\item For two
contractive maps $f:X\To Y$ and $g:Y\To Z$ we have
$\overline{gf}=\bar{g}\bar{f}$.\end{enumerate}
\end{teo}

\begin{dem}
For each element $x\in\bar{X}$, choose an expression of $x$ as a
convex combination of elements of $X$, like
$x=\sum_{i}a_{i}x_{i}$. If we want $\bar{f}$ to be contractive it
must be defined like $\bar{f}(x)=\sum_{i}a_{i}f(x_{i})\in\bar{Y}$.
This proves uniqueness. For existence we must check that, so
defined, $\bar{f}$ is contractive. We take $x,y\in\bar{X}$, and
their corresponding expressions $x=\sum
    a_{i}x_{i}$ and $y=\sum b_{j}y_{j}$ with $x_{i},y_{j}\in X$:
    $$ d(\bar{f}(x),\bar{f}(y)) = \bigoplus a_{i}b_{j}d(f(x_{i}),f(y_{j}))
    \leq \bigoplus a_{i}b_{j}d(x_{i},y_{j}) = d(x,y) $$
    If $f$ is immersion then the inequality turns into an equality, and we deduce that
    $\bar{f}$ is an immersion. Property (2) is trivial
    and from this, using $f^{-1}$, we deduce that if $f$ is isometry so is $\bar{f}$.
\end{dem}

As a corollary, we get that the convex closure of a metric space is unique, up to 
isometry, since
if $X\subseteq X_{1},X_{2}$ are two convex closures of $X$, then
$1_{X}$ extend to an isometry $f:X_{1}\To X_{2}$.\\

In the sequel $conv(X)$ will denote a convex closure of $X$. We finish by stating some elementary properties of convex spaces and
convex closures.\\

Let $X$ and $Y$ be convex metric spaces over $B$ and $U\subseteq X$. Then, 
the following hold:
\begin{enumerate}
\item  The set of all convex combinations of elements of $U$ in $X$
is a convex closure of $U$(In this situation, the notation
$conv(U)$ will refer to this set). \item If $f:X\longrightarrow Y$
is contractive, then $f(conv(U))=conv(f(U))$.
\item If $X_{1},\ldots,X_{n}$ are metric spaces over $B$, then $conv(X_{1})\times\cdots\times conv(X_{n})$ is a
convex closure of $X_{1}\times\cdots\times X_{n}$.\end{enumerate}

\section{CFG-spaces}

\begin{defi}
A metric space $X$ over $B$ will be said to be a CFG-space
(convex finitely generated space) if it is the convex closure of a
finite subspace.
\end{defi}

Observe that
\begin{enumerate}
\item If $X$ is a CFG-space and $f:X\longrightarrow Y$ is contractive, then $f(X)$ is a
CFG-space.
\item The product of a finite number of CFG-spaces is a CFG-space.
\end{enumerate}

For technical reasons, it is convenient to work with pointed
metric spaces. $(X,0)$ is said to be a (pointed) metric space if
$X$ is a metric space over $B$ and $0\in X$. A map
$f:(X,0)\longrightarrow (Y,0')$ will mean a map
$f:X\longrightarrow Y$ such that $f(0)=0'$, and expressions like
$x\in (X,0)$ will mean simply $x\in X$.\\

Throughout this section, we fix a convex metric space $(X,0)$. By
Theorem~\ref{inclusion en un modulo}, it is not restrictive to
suppose that $X$ is a metrizable convex subset of a $B$-module $M$
and that 0 is the zero element of $M$. In $(X,0)$ we will use the
following notations:

\begin{itemize}
\item For $x\in X$, $|x|:=d(0,x)$.
\item If $x_{1},\ldots,x_{n}\in X$ and $a_{1},\ldots,a_{n}\in B$
are such that $a_{i}a_{j}=0$ whenever $i\neq j$, then we have an
element of $X$: $$a_{1}x_{1}+\cdots+a_{n}x_{n}=
a_{0}0+a_{1}x_{1}+\cdots+a_{n}x_{n}$$ where $a_{0}=1+a_{1}+\cdots
a_{n}$ (note that the right expression represents an element of
$X$ since $a_{0}\oplus\cdots\oplus a_{n}=1$ and $X$ is a convex
space). Such a combination will be called an orthogonal
combination. As a particular case, $ax =
ax+\bar{a}0$ for $x\in X$, $a\in B$.
\item For $x\in X$, $Bx := \{ax : a\in B\} = conv(0,x)$.
\item For $x,y\in X$, $x\star y := \overline{d(x,y)}x$.
\end{itemize}

Note that any contractive map $f:(X,0)\longrightarrow (Y,0')$
preserves orthogonal combinations. In the following lemma, we state 
some elementary properties:

\begin{lema}\label{propiedades de los espacios centrados}
Let $x,y\in X$ and $a,b\in B$. Then:
\begin{enumerate}
\item The maps $||:(X,0)\longrightarrow (B,0)$ and
$x\star\_:(X,0)\longrightarrow (X,0)$ are contractive, so both
preserve orthogonal combinations.
\item $ax=bx$ if and only if $a +b \in \overline{|x|}B$ (if and only if $a+\overline{|x|}B = b+\overline{|x|}B$).
\item $ax=0$ if and only if $a\leq\overline{|x|}$, and $ax=x$ if and only if $a\geq |x|$.
\item The operation $(\star)$ is commutative.
\end{enumerate}
\end{lema}

\begin{dem}
For property 1, the function $x\star\_$ can be expressed as a
composition of $y\mapsto d(x,y)$, $b\mapsto \bar{b}$ and $b\mapsto
bx$ and all of them are contractive.

 Property 2: Suppose $X$ is a metrizable subset of a $B$-module. Then,
 $ax=bx$ if and only if $a+b\in Ann(x-0) = \overline{d(x,0)}B$.

 Property 3 follows from 2.

 For property 4, $x\star y = y\star x$ if and only if
$\overline{d(x,y)}x+d(x,y)0 = \overline{d(x,y)}y + d(x,y)0$. This
equality is easily checked verifying that the distance between the
two terms is zero using Proposition~\ref{convexidad de la
distancia}.
\end{dem}

\begin{lema}\label{significado de estrella}
$Bx\cap By = B(x\star y)$ for all $x,y\in X$.
\end{lema}

\begin{dem}
Just by the definition of $\star$ we have $B(x\star y)\subseteq
Bx$, and symmetrically, since $\star$ is commutative, $B(x\star
y)\subseteq By$, so one inclusion is proved. Now suppose $u\in
Bx\cap By$. Then $u=ax=by$, and if we call $c=ab$ then
$cx=bax=bu=bby=u=aax=au=aby=cy$. Thus, $cx=u=cy$, and that implies
$c\in Ann(x-y) = \overline{d(x,y)}B$ (suppose $X$ is a subset of a
module) and $u = cx = c\overline{d(x,y)}x = c(x\star y)$.
\end{dem}

\begin{prop}
For two elements $x,y\in X$ the following are equivalent:
\begin{enumerate}
\item $x\star y = 0$
\item $Bx\cap By = \{0\}$
\item $d(x,y) = |x|\vee |y|$
\end{enumerate}
In this case, $x$ and $y$ will be said to be orthogonal and we
will write $x\perp y$.
\end{prop}

\begin{dem}
($1\Leftrightarrow 2$) is a direct consequence of
Lemma~\ref{significado de estrella}. For ($1\Leftrightarrow 3$),
we have (1) if and only if $0= |x\star y| =\overline{d(x,y)}|x|$
and, by symmetry, if and only if
$\overline{d(x,y)}|x|=0=\overline{d(x,y)}|y|$, which is equivalent
to $|x|,|y|\leq d(x,y)$, and $|x|\vee |y|\leq d(x,y)$. The
converse of the latter inequality is always true by axiom 3 of
Boolean metric spaces.
\end{dem}

For $x,y\in (B,0)$, we have $x\star y =
\overline{d(x,y)}x=(x+y+1)x = xy$, so this concept of
orthogonality corresponds to \emph{disjointness} in $B$.

\begin{defi}
A finite subset $R\subseteq X$ will be said to be orthogonal if
every two different elements in $R$ are orthogonal, and $0\not\in
R$. If, moreover, $X=conv(R\cup\{0\})$, $R$ will be said to be a
reference system or a referential of $(X,0)$.
\end{defi}

\begin{prop}\label{coordenadas}
Let $R = \{x_{1},\ldots,x_{n}\}$ be a reference system of $(X,0)$
and $x\in X$. There is a unique tuple $(a_{1},\ldots,a_{n})\in
B^{n}$ satisfying the three following properties:
\begin{enumerate}
\item $a_{i}a_{j}=0$ whenever $i\neq j$.
\item $\sum_{1}^{n}a_{i}x_{i}=x$.
\item $a_{i}\leq |x_{i}|$ for $i=1,\ldots,n$.
\end{enumerate}
Such a tuple will be called the tuple of coordinates of $x$ with
respect to $R$.
\end{prop}

\begin{dem}
Uniqueness: If $\sum_{1}^{n}a_{i}x_{i} = \sum_{1}^{n}b_{i}x_{i}$
in those conditions, multiplying by $a_{i}b_{j}$, $i\neq j$, we
obtain $a_{i}b_{j}x_{i} = a_{i}b_{j}x_{j}\in Bx_{i}\cap
Bx_{j}=\{0\}$ so for each $i$, $a_{i}x_{i} =
a_{i}(\sum_{j}b_{j})x_{i} = a_{i}b_{i}x_{i}$, and symmetrically
$b_{i}x_{i} = a_{i}b_{i}x_{i} = a_{i}x_{i}$, so by
Lemma~\ref{propiedades de los espacios centrados} $a_{i}+b_{i}\in
\overline{|x_{i}|}B$, and also $a_{i}+b_{i}\in|x_{i}|B$ since the
$a_{i}$ and $b_{i}$'s are assumed to verify property 3. So
$a_{i}+b_{i}=0$ for all $i$.

Existence: Since $X=conv\{0,x_{1}\ldots,x_{n}\}$, we can find
$b_{1},\ldots,b_{n}\in B$ verifying 1 and 2. Now set
$a_{i}=|x_{i}|b_{i}$. The $a_{i}$'s satisfy trivially 1 and 3.
Using Lemma~\ref{propiedades de los espacios centrados} we deduce
from $a_{i}+b_{i} = \overline{|x_{i}|}b_{i}\in\overline{|x_{i}|}B$
that $a_{i}x_{i}=b_{i}x_{i}$ for all $i$. So
$\sum_{1}^{n}a_{i}x_{i} = \sum_{1}^{n}b_{i}x_{i}=x$.
\end{dem}

\begin{prop}\label{definicion de funciones sobre referenciales}
Let $R = \{x_{1},\ldots,x_{n}\}$ be a referential of $(X,0)$
and $(Y,0')$ a convex metric space. Then, $f:R\longrightarrow Y$
is extensible to a (unique) contractive map
$\hat{f}:(X,0)\longrightarrow (Y,0')$ if and only if
$|f(x_{i})|\leq |x_{i}|$ for $i=1,\ldots,n$.
\end{prop}

\begin{dem}
Define $f$ on $R\cup\{0\}$ by $f(0)=0'$. By Theorem \ref{extension
a la envoltura convexa}, $f$ admits such an extension if and only
if it is contractive. If $f$ is contractive, it is clear that
$|f(x_{i})|\leq |x_{i}|$ for $i=1,\ldots,n$, so one way is proved.
Conversely, suppose $|f(x_{i})|\leq |x_{i}|$ for every $i$. Then,
for all $i\neq j$, since $x_{i}$ and $x_{j}$ are orthogonal, we
have last equality in $d(f(x_{i}),f(x_{j}))\leq |f(x_{i})|\vee
|f(x_{j})|\leq |x_{i}|\vee |x_{j}| = d(x_{i},x_{j})$.
\end{dem}

We check now that any CFG-space possesses a reference system.

\begin{teo}\label{construccion de referenciales}
    Suppose $X=conv\{0,x_{1},\dots,x_{n}\}$ and that the set
    $\{x_{1},\dots,x_{s}\}$ is orthogonal. Then, there exist
    $a_{s+1},\dots,a_{n}\in B$ such that
    $\{x_{1},\dots,x_{s},a_{s+1}x_{s+1},\dots,a_{n}x_{n}\}\setminus\{0\}$
    is a referential of $(X,0)$.
\end{teo}

\begin{dem}
    Let $r=card\{(i,j) : x_{i}\star x_{j}\neq 0 \}$. We make
    induction on $r$. We suppose that the theorem holds for any value lower than $r>0$.
    We take $x_{i},x_{j}$ with
    $x_{i}\star x_{j} \neq 0$ and suppose, without loss of generality
    that $i,s<j$. Let $a:=d(x_{i},x_{j})$. Since
    $ax_{j}\star x_{i}=a\overline{d(x_{i},x_{j})}x_{i}=0$, we have
    $ax_{j}\perp x_{i}$. Also, $x_{j} =
    a(ax_{j}) + \bar{a}x_{i}$, and from this we deduce that $conv\{0,x_{i},x_{j}\}
    = conv\{0,x_{i},ax_{j}\}$ and therefore:
    $$X =
    conv\{0,x_{1},\dots,x_{j-1},ax_{j},x_{j+1},\dots,x_{n}\}$$
    Making use of the induction hypothesis, the proof is complete (in
    this system of generators there is at least one orthogonal pair
    more, since $x_{i}\perp ax_{j}$).
\end{dem}

\begin{coro}
    Let $\{x_{1},\dots,x_{s}\}$ be an orthogonal subset of the CFG-space $(X,0)$.
    Then, there exist $x_{s+1},\dots,x_{n}\in X$ such that
    $\{x_{1},\dots,x_{n}\}$ is a referential of $(X,0)$.
\end{coro}

In particular, any CFG-space $(X,0)$ has a reference system.

\begin{defi}For $U\subseteq X$, $U^{\perp} = \{x\in X : x\perp y
\ \forall y\in U\}$.\end{defi}

\begin{prop}\label{ortogonal de CFG es CFG}
For two CFG-spaces $(U,0)\subseteq (X,0)$, the space $U^{\perp}$
is a CFG-space and $conv(U\cup U^{\perp}) =X$.
\end{prop}

\begin{dem}
    Let $\{x_{1},\dots,x_{m}\}$ be a reference system of $(U,0)$
    that we can extend to a reference system $\{x_{1},\dots,x_{n}\}$
    of $(X,0)$. We prove that
     $U^{\perp} = conv\{0,x_{m+1},\dots,x_{n}\}$.
     Take $x\in U^{\perp}$, $x=\sum_{1}^{n}a_{i}x_{i}$.
    Then, for $j=1,\ldots,m$ we have $0=x_{j}\star x =
    \sum_{1}^{n}a_{i}(x_{i}\star x_{j}) = a_{j}x_{j}$. Hence,
    $x=\sum_{1}^{m}a_{i}x_{i}$.
\end{dem}

\begin{prop}\label{suma ortogonal de aplicaciones}
    Let $(U,0)\subset (X,0)$ be two CFG-spaces, $(Y,0')$
    a convex metric space and $f:(U,0)\rightarrow (Y,0')$ and
    $g:(U^{\perp},0)\rightarrow (Y,0')$
    contractive maps. Then, there is a unique contractive map $f\perp
    g:(X,0)\rightarrow (Y,0')$ that extends $f$ and $g$.
\end{prop}

\begin{dem}
    We take a referential of $(U,0)$ and another one of $(U^{\perp},0)$. The
    union is a referential of $(X,0)$.
    Applying Theorem~\ref{definicion de funciones sobre referenciales}
    the proposition is proved.
\end{dem}

\section{Algebraic Geometry over CFG-rings}

\begin{defi}\label{CFG-anillo}
A regular ring $A$ is said to be a CFG-ring if, equipped with its
modular metric, it is a CFG-space over $B(A)$.
\end{defi}

In this case, $A^{n}$ (which is a metrizable $A$-module for which
the product metric and the modular metric coincide) is a CFG-space
over $B(A)$ too. If $p$ is a prime number, any $p$-ring is a
CFG-ring because if $A$ is a $p$-ring, then
$A=conv\{0,1,\ldots,p-1\}$. A proof of this fact can be found
in~\cite{Zem} (Corollary 1). There are CFG-rings that are not
$p$-rings. For instance, take $K$ a finite field and $\Omega$ a
set. Then, $K^{\Omega}$ is regular and it is easy to see that the
set of constant tuples constitute a finite system of
generators of $K^{\Omega}$, so $K^{\Omega}$ is a CFG-ring.
The aim of this section is to prove Theorem~\ref{variedades algebraicas}.

\begin{lema}\label{polinomica implica convexa}
Let $R$ be a ring and $f:R^{n}\To R$ a polynomial function. For
every $x_{1},\ldots,x_{m}\in R^{n}$ and every $e_{1},\ldots,e_{n}\in B(R)$ such that
$e_1\oplus\cdots\oplus\cdots e_n=1$,
 we have
$f(\sum_{i}e_{i}x_{i}) = \sum_{i}e_{i}f(x_{i})$.
\end{lema}

\begin{dem}
Let $S$ be the set of all maps $g:R^{n}\To R$ verifying the
conclusion of the lemma. It is straightforward to check that the
projections $\pi_{i}:R^{n}\To R$ are in $S$ (that proves the lemma
for the polynomials $X_{1},\ldots,X_{n}$), that constant maps are
in $S$, and that the sums and products of maps in $S$ lie in $S$.
Any polynomial map is a sum of products of constants and the
variables $X_{i}$'s.
\end{dem}

\begin{lema}\label{contractiva sii polinomica en A}
Let $A$ be a CFG-ring. A map $f:A^{n}\To A^{m}$ is contractive if
and only if it is a polynomial map.
\end{lema}

\begin{dem}
    $f$ is contractive if and only if all its components are, and
    the same holds about $f$ being a polynomial map, so we can
    assume $m=1$. The `if' part is a consequence of
    Lemma~\ref{polinomica implica convexa}, Theorem~\ref{convexa sii
    contractiva}, and Proposition~\ref{combinacion convexa=combinacion
    lineal}. For the `only if' part, we assume that $f(0)=0$ (it is 
    sufficient
    to prove this case, just considering for an
    arbitrary $f$, the composition $h\circ f$ where $h:A\To A$ is
    $h(x) = x+f(0)$). Take $\{x_{1},\ldots,x_{r}\}$ a referential of
    $(A^{n},0)$. We prove first the case $n=1$:

    \emph{Case $n=1$}. Consider the polynomial
    $g_{i}(x) = x\prod_{j\neq i}(x-x_{j})$  for $i=1,\ldots,r$. Then,
    \begin{eqnarray*}e(g_{i}(x_{i})) &=& e(x_{i})\prod_{j\neq i}e(x_{i}-x_{j}) =
    e(x_{i})\prod_{j\neq i}d(x_{i},x_{j})\\ &=& e(x_{i})\prod_{j\neq
    i}\left(|x_{i}|\vee|x_{j}|\right)
     = |x_{i}| = e(x_{i})\end{eqnarray*} Therefore, for each $i=1,\ldots,n$ there is a
    unit $a_{i}$ of $A$ such that $g_{i}(x_{i}) =
    a_{i} e(x_{i})$. Consider the polynomial map $g:A\To A$
    given by $$g(x) = \sum_{i=1}^{r}a_{i}^{-1}f(x_{i})g_{i}(x)$$
    Since $A=conv\{0,x_{1},\ldots,x_{r}\}$ and $f$ and $g$ are contractive, we prove that $g(x)$ and $f(x)$ coincide for
    $x=0,x_{1},\ldots,x_{r}$. It is clear that
    $g(0) = 0 = f(0)$ because $g_{i}(0) = 0$. For $x=x_{j}$, $$g(x_{j}) =
    \sum_{i=1}^{r}a_{i}^{-1}f(x_{i})g_{i}(x_{j})$$ and since $g_{i}(x_{j})=0$ whenever
    $i\neq j$, $$g_{i}(x_{j})=a_{j}^{-1}f(x_{j})g_{j}(x_{j}) = f(x_{j})e(x_{j})$$ and this equals
    $f(x_{j})$ because $|f(x_{j})|\leq
    |x_{j}|= e(x_{j})$.

    \emph{General case}: As a consequence of the case $n=1$, we
    find that $e:A\To A$ is a polynomial map, and therefore,
    for $v\in A^{n}$, the map $d(\sim,v):A^{n}\To A$ is
    polynomial too, since if $v=(a_{1},\ldots,a_{n})$ then
    $d(x,v) = e(x_{1}-a_{1})\vee \cdots \vee
    e(x_{n}-a_{n})$ (recall that $x\vee y = x + y - xy$
    for $x,y\in B(A)$). Hence, we can construct polynomial maps
    for $i=1,\ldots,r$ given by $G_{i}(x) = |x|\prod_{j\neq
    i}d(x,x_{j})$. We define $G(x) = \sum_{i=1}^{r}f(x_{i})G_{i}(x)$.
    and we are going to see that $G$ and $f$ coincide on
    $\{0,x_{1},\ldots,x_{r}\}$, so that, since both are
    contractive and this set generates $A^{n}$, that will prove
    that $f=G$, so $f$ is a polynomial map. It is clear that $G(0) = 0
    = f(0)$ because $G_{i}(0)=0$ for all $i$. Since $G_{i}(x_{j}) = 0$
    if $i\neq j$, \begin{eqnarray*}G(x_{i}) &=& f(x_{i})G_{i}(x_{i}) =
    f(x_{i})|x_{i}|\prod_{i\neq j}d(x_{i},x_{j})\\ &=& f(x_{i})|x_{i}|\prod_{j\neq i}|x_{i}|\vee|x_{j}|
    = f(x_{i})|x_{i}| = f(x_{i})\end{eqnarray*} where the last equality follows from
    the fact that $|f(x_{i})|\leq |x_{i}|$.
\end{dem}

\begin{lema} \label{conjunto de soluciones}
    Let $(X,0) = conv(H)$ be a convex metric space with $0\in
    H$,
    and let $f:(X,0)\rightarrow (Y,0')$ be contractive. Then
    $$ f^{-1}(0') = conv\{0, \overline{|f(x)|}x : x\in H\}$$
\end{lema}
\begin{dem}
    One inclusion is clear because all the elements that
    appear in the right term are in the convex set $f^{-1}(0')$.
    For the converse, if $x\in f^{-1}(0')$, in particular
    it is in $X$, so it can be expressed like
    $x=\sum_{i=1}^{n}a_{i}x_{i}$ with $x_{i}\in H$,
    and $a_{i}a_{j}=0$ whenever $i\neq j$.
    Then, $$0 = |f(x)| =
    a_{1}|f(x_{1})| \oplus \cdots \oplus a_{n}|f(x_{n})|,$$ and
    therefore
    $a_{i}|f(x_{i})| = 0$,  so $a_{i}=a_{i}\overline{|f(x_{i})|}$
    for all $i=1,\dots,n$. Finally,
    $$x=\sum_{i=1}^{n}a_{i}x_{i}=\sum_{i=1}^{n}a_{i}\overline{|f(x_{i})|}x_{i}\in
    conv\{0, \overline{|f(x_{1})|}x_{1}, \dots,
    \overline{|f(x_{n})|}x_{n}\}.$$
\end{dem}

Note that a convex subset of a CFG-space need not be a CFG-space. For
instance, $B=conv\{0,1\}$ is a CFG-space, and those ideals of $B$
that are not finitely generated are convex subsets that are not
CFG-spaces.

\begin{lema}\label{CFG es nucleo}
    Let $X$ be a CFG-space. Then $Y\subseteq X$ is a CFG-space if and only if
    there exists a contractive map $f:X\rightarrow B$ such that $Y =
    f^{-1}(0)$.
\end{lema}

\begin{dem}
    One way is a direct consequence of Lemma~\ref{conjunto de soluciones}.
    For the converse, suppose $Y$ is a CFG-space. If $Y=\emptyset$ it is trivial and
    if not, take $0'\in Y$ and
    $\{u_{1},\dots,u_{k}\}$ a referential in $(Y,0')$
    that we extend to a referential of $(X,0')$,
    $\{u_{1},\dots,u_{n}\}$. By Theorem~\ref{definicion de funciones sobre
    referenciales} we can define a contractive map
    $f:(X,0')\rightarrow (B,0)$ such that $f(u_{i}) = 0$ if $i \leq k$ and $f(u_{i}) =
    |u_{i}|$ if $i>k$. It is clear that $ Y \subset f^{-1}(0)$ and
    for the other inclusion suppose $x\in f^{-1}(0)$ has
    coordinates $(a_{1},\dots,a_{n})$. Then
    $$0 = f(x) = f( \sum_{i=1}^{n}a_{i}u_{i}) = \bigoplus_{i=1}^{n}a_{i}f(u_{i}) =
    \bigoplus_{i=k+1}^{n}a_{i}|u_{i}| = \bigoplus_{i=k+1}^{n}a_{i}$$
    so $a_{i} = 0$ for $i>k$, and $x=\sum_{i=1}^{k}a_{i}u_{i}\in
    conv\{0',u_{1},\dots,u_{k}\} = Y$.
\end{dem}

\begin{coro}\label{interseccion de CFGs}
    If $Y,Z$ are CFG-spaces contained in the space $X$, then $Y\cap Z$
    is a CFG-space.
\end{coro}

\begin{dem}
    If $Y = f^{-1}(0)$ and $Z = g^{-1}(0)$ with $f,g:conv(Y\cup Z)\rightarrow B$
    contractive maps, then $Y\cap Z = (f\vee g)^{-1}(0)$.
\end{dem}

\begin{coro}\label{preimagen de un CFG}
    Let $f:X\rightarrow Y$ be a contractive map between CFG-spaces.
    If $Z\subset Y$ is a CFG-space, then $f^{-1}(Z)$ is a CFG-space too.
\end{coro}

\begin{dem}
    Let $g:Y\rightarrow B$ be such that $K = g^{-1}(0)$. Then
    $f^{-1}(K) = (g\circ f)^{-1}(0)$, so it is a CFG-space.
\end{dem}

\begin{teo}\label{variedades algebraicas}
Let $A$ be a CFG-ring.
\begin{enumerate}
\item A subset $U\subseteq A^{n}$ is an algebraic variety
if and only if $U$ is a CFG-metric subspace of $A^{n}$.

\item A map $f:U\To V$ between two algebraic varieties is a polynomial map if and only if it is
contractive.
\end{enumerate}
\end{teo}

\begin{dem}
If $U$ is an algebraic variety then,
$U=\bigcap_{1}^{k}f_{k}^{-1}(0)$ where $f_{i}:A^{n}\To A$ are
polynomial maps, and therefore, by Lemma~\ref{contractiva sii
polinomica en A}, contractive maps. Using Lemma~\ref{CFG es
nucleo} and Corollary~\ref{interseccion de CFGs} we deduce that
$U$ is a CFG-space. Conversely, If $U$ is a CFG-space, by
Lemma~\ref{CFG es nucleo}, there is a contractive map
$f:A^{n}\To B(A)$ with $U=f^{-1}(0)$. Then $U=g^{-1}(0)$ where $g$
is the composition $A^{n}\To B(A)\hookrightarrow A$, that is
contractive and therefore polynomial, again by
Lemma~\ref{contractiva sii polinomica en A}.

Suppose $f:U\To V$ is contractive, choose some $u\in U$ and consider
$f:(U,u)\To (V,f(u))$ and $k:(U^{\perp},u)\To
(V,f(u))$ the constant map. Since $U$ and $A^{n}$ are CFG-spaces and $V$ is convex,
we can consider, by Proposition~\ref{suma ortogonal de
aplicaciones}, $f\perp k:A^{n}\To A^{m}$ that is contractive, and
therefore, a polynomial map, that extends $f$. The converse is a
direct consequence of Lemma~\ref{contractiva sii polinomica en
A}.
\end{dem}

\section{Structure Theorem for CFG-spaces}

We shall classify now CFG-spaces up to isomorphism. Reference
systems do not give good isomorphism invariants, since they are
not \emph{unique up to isometry} (for instance, $\{1\}$ and
$\{a,\bar{a}\}$ are non-isometric referentials of $(B,0)$). The
right concept for this purpose is the following:

\begin{defi}
A referential $\{x_{1},\ldots,x_{n}\}$ of $(X,0)$ is said to be a
base of $(X,0)$ if $|x_{1}|\geq |x_{2}|\geq\cdots\geq |x_{n}|$.
\end{defi}

We will prove that there exists a base for any pointed CFG-space
$(X,0)$, and that they are unique in the sense of
Theorem~\ref{unicidad de las bases} below. We prove uniqueness
first, and existence afterwards.

\begin{defi}
    Let $k>0$ be an integer and $X$ a metric space over $B$. The k-ideal of
    $X$ (denoted by $I_{k}(X)$) is the ideal of $B$ generated by
    $$\{\prod_{0\leq i < j \leq k}d(u_{i},u_{j}) : u_{0},\dots,u_{k}\in X\}.$$
    If $I_{k}(X)$ is principal, we will denote by $\alpha_{k}(X)$
    its generator.
\end{defi}

\begin{lema}\label{amplitudes de la envoltura}
    If $X = conv(H)$, then $I_{k}(H) = I_{k}(X)$ for all $k\in
    \mathbf{N}$.
\end{lema}

\begin{dem}
    If $U$ is any Boolean metric space, the map $f_{U}:U^{k+1}\rightarrow B$ given by
    $f_{U}(u_{0},\dots,u_{k}) = \prod_{0\leq i < j \leq k}d(u_{i},u_{j})$ is
    contractive because it is a composition of \emph{distance functions}
    and a polynomial function (the product). With this notation,
    $I_{k}(U)$ is the ideal generated by the image of $f_{U}$, and
    \begin{eqnarray*} Im(f_{X}) &=& f_{X}(X^{k+1}) =
    f_{X}(conv(H^{k+1}))=conv(f_{X}(H^{k+1}))\\ &=&
    conv(Im(f_{H}))\end{eqnarray*}
    so both images generate the same ideal.
\end{dem}

\begin{lema}\label{existencia de las amplitudes}
Let $X$ be a CFG-space. Then $I_{k}(X)$ is principal for all $k\in
\mathbf{N}$ and there exists $n\in\mathbf{N}$ such that
$I_{k}(X)=0$ for all $k\geq n$. Hence, $\alpha_{k}(X)$ exists for
all $k\in \mathbf{N}$ and $\alpha_{k}(X)=0$ for $k\geq
n$.\end{lema}

\begin{dem}
Suppose $X=conv(H)$ with $H$ finite. Then, $I_{k}(X) = I_{k}(H)$
is always a finitely generated ideal of $B$, so it is principal,
and if we take $n=card(H)$, $0=I_{k}(H)=I_{k}(X)$ if $k\geq n$.
\end{dem}

\begin{lema}\label{amplitudes de una base}
    Let $\{x_{1},\dots,x_{n}\}$ be a base of $(X,0)$. Then,
    $\alpha_{k}(X) = |x_{k}|$ for $k \leq n$ and $\alpha_{k}(X) = 0$
    if $k>n$.
\end{lema}

\begin{dem}
     By Lemma~\ref{amplitudes de la envoltura}, $\alpha_k(X)=\alpha_{k}(H)$ where
     $H=\{0,x_{1},\dots,x_{n}\}$. If $k>n$,
    it is trivial that $\alpha_{k}(H) = 0$. If $k\leq n$, call $y_{i}$'s to the reordering
    of the $x_{i}$'s such that $0=|y_{0}| \leq |y_{1}| \leq \cdots \leq |y_{n}|$
    ($y_{r}=x_{n-r+1}$ if $r>0$). For $i<j$ we have
    $d(y_{i},y_{j}) = |y_{i}|\vee |y_{j}| = |y_{j}|$, by orthogonality.
    We wonder whether $I_{k}(H) = |x_{k}|B (=|y_{n-k+1}|B)$.
    One
    inclusion is because $$|y_{n-k+1}| = \prod_{n\geq i> n-k}|y_{i}|=
    \prod_{n\geq i>j\geq n-k}d(y_{i},y_{j})$$ is one of the generators of
    $I_{k}(H)$. For the other inclusion we shall check that all the
    generators of $I_{k}(H)$ are in the ideal $|y_{n-k+1}|B$. Take
    $U=\{u_{0},\ldots,u_{k}\}\subseteq H$. By a cardinality argument,
    there must exist indices $r<s\leq n-k+1$ such that
    $y_{r},y_{s}\in U$, so $$\prod_{0\leq i < j \leq
    k}d(u_{i},u_{j})\leq d(y_{r},y_{s}) = |y_{s}|\leq
    |y_{n-k+1}|$$
\end{dem}

\begin{teo}\label{unicidad de las bases}
If $\{x_{1},\ldots,x_{n}\}$ is a base of $(X,0)$ and
$\{y_{1},\ldots,y_{m}\}$ is a base of $(X,0')$, then $n=m$ and
$|x_{i}|=|y_{i}|$ for $i=1,\ldots,n$. Moreover, there exists an
isometry $f:(X,0)\To (X,0')$ such that $f(x_{i}) = y_{i}$ for
$i=1,\ldots,n$.
\end{teo}

\begin{dem}
By Lemma~\ref{amplitudes de una base}, we know that
$n=\max\{k:\alpha_{k}(X)\neq 0\} = m$ and $|x_{i}| = \alpha_{i}(X)
= |y_{i}|$. About the last assertion, there exists a contractive
map $f:(X,0)\To (X,0')$ such that $f(x_{i}) = y_{i}$ for
$i=1,\ldots,n$, by virtue of Proposition~\ref{definicion de
funciones sobre referenciales}. It is an isometry because we can
find its inverse in an analogue way $g:(X,0')\To (X,0)$ with
$g(y_{i})=x_{i}$.
\end{dem}

\begin{lema}
    If $(V,0)\subset (X,0)$ are CFG-spaces and $V^{\perp} =
    \{0\}$, then $V=X$.\label{cancelacion restringida}
\end{lema}

\begin{dem}
    A reference system of $(V,0)$, $\{x_{1},\dots,x_{m}\}$ can be extended to a
    referential of $(X,0)$, $\{x_{1},\dots,x_{n}\}$. Then,
    $x_{m+1},\dots,x_{n}\in V^{\perp} = \{0\}$, so
    $V=conv\{0,x_{1},\dots,x_{m}\}=conv\{0,x_{1},\dots,x_{n}\}=X$.
\end{dem}

\begin{lema}\label{Weierstrass}
Let $X$ be a CFG-space and $f:X\To B$ contractive. Then, there
exists $u\in X$ such that $f(u)=\max\{f(x):x\in X\}$.
\end{lema}

\begin{dem}
Suppose $X=conv\{x_{0},\ldots,x_{n}\}$. The set $f(X)\subseteq B$
is closed under the operation $(\vee)$, because it is convex and
$a\vee b = aa + \bar{a}b$. Hence, there exists $u\in X$ such that
$f(u)=f(x_{0})\vee\cdots\vee f(x_{n})$. If $x\in X$, we express it
as a convex combination $x=\sum_{i}a_{i}x_{i}$ and
$f(x)=\bigoplus_{i}a_{i}f(x_{i}) = \bigvee_{i}a_{i}f(x_{i})\leq
\bigvee_{i}f(x_{i}) = f(u)$.
\end{dem}

\begin{teo}\label{existencia de bases}
Any CFG-space $(X,0)$ possesses a base.
\end{teo}

\begin{dem}
    We define by recursion a sequence $(x_{n})_{n=1}^{\infty}$ in $X$
    and a sequence $(U_{n})_{n=1}^{\infty}$ of CFG-spaces
    contained in $X$:
    \begin{itemize}
        \item $x_{1}$ is such that $|x_{1}| = \max\{|x| : x\in X\}$; $U_{1}=conv\{0,x_{1}\}$
        \item Given $x_{i}$ and $U_{i}$ for $i<n$, we take $x_{n}$ such that $|x_{n}| = \max\{|x| : x\in U_{n-1}^{\perp}\}$
        and $U_{n}:=conv\{0,x_{1},\dots,x_{n}\}$
    \end{itemize}
    Note that those maximums exist by virtue of Lemma~\ref{Weierstrass}, since
    $U_{n-1}^{\perp}$ is a CFG-space by Proposition~\ref{ortogonal de CFG es CFG}.
    The $x_{i}$'s form an orthogonal set and verify $|x_{i}| \geq |x_{j}|$ whenever $i<j$.
    Therefore, $\{x_{1},\ldots,x_{n}\}\setminus\{0\}$ is a
    base of $(U_{n},0)$. Since $X$ is a CFG-space, by Lemma~\ref{existencia de las
    amplitudes},
    there must exist some $k>0$
    with $0 = \alpha_{k}(X) \geq \alpha_{k}(U_{k}) = |x_{k}|$. So, taking $r$
    the largest integer such that $|x_{r}|\neq 0$, we have, just by
    the definition of $x_{r+1}=0$, that
    $U_{r}^{\perp} = \{0\}$.
    Therefore, by Lemma~\ref{cancelacion restringida},
    $U_{r} = X$ and we have already shown
    that
    $\{x_{1},\dots,x_{r}\}\setminus\{0\}$ is a base of $(U_{r},0)$.
\end{dem}

\begin{teo}\label{isometricos=mismas
amplitudes}
    Two CFG-spaces $X$ and $Y$ are isometric if and only if $\alpha_{k}(X) =
    \alpha_{k}(Y)$ for all $k\in\mathbf{N}$.
\end{teo}

\begin{dem}
    Suppose $\alpha_{k}(X) =
    \alpha_{k}(Y)$ for all $k$. Choose $0\in
    X$, $0'\in Y$ and bases
    $\{x_{1},\dots,x_{n}\}$ and $\{y_{1},\dots,y_{m}\}$
    of $(X,0)$ and $(Y,0')$ respectively. Then,
    $n=\max\{k:\alpha_{k}(X)=\alpha_{k}(Y)\neq 0\} = m$ and
    we can construct an
    isometry like in the proof of Theorem~\ref{unicidad de las bases}.
\end{dem}

\end{document}